\documentclass[a4paper,12pt,reqno]{amsart}
\usepackage{amssymb}
\usepackage{graphicx}
\usepackage{amsmath}
\usepackage{color}
\usepackage{enumitem}
   \nonstopmode \numberwithin{equation}{section}
\newtheorem{thm}{Theorem}[section]
 \newtheorem{cor}[thm]{Corollary}
 \newtheorem{lem}[thm]{Lemma}
 \newtheorem{prop}[thm]{Proposition}
 \theoremstyle{definition}
 \newtheorem{defn}[thm]{Definition}
 \theoremstyle{remark}
 \newtheorem{rem}[thm]{Remark}
 \newtheorem{rems}[thm]{Remarks}
  \newtheorem{ques}[thm]{Question}
 \newtheorem{ex}[thm]{Example}

 \newtheorem{prob}[thm]{Problem}
 \numberwithin{equation}{section}

\def\be{\begin{equation}}
\def\ee{\end{equation}}

\newcommand{\blem}{\begin{lem}}
\newcommand{\elem}{\end{lem}}
\newcommand{\bthm}{\begin{thm}}
\newcommand{\ethm}{\end{thm}}
\newcommand{\bcor}{\begin{cor}}
\newcommand{\ecor}{\end{cor}}
\newcommand{\beg}{\begin{ex}}
\newcommand{\eeg}{\end{ex}}
\newcommand{\bprop}{\begin{prop}}
\newcommand{\eprop}{\end{prop}}
\newcommand{\bdefn}{\begin{defn}}
\newcommand{\edefn}{\end{defn}}
\newcommand{\bprob}{\begin{prob}}
\newcommand{\eprob}{\end{prob}}
\newcommand{\bpf}{\begin{proof}}
\newcommand{\epf}{\end{proof}}

\newcommand{\brem}{\begin{rem}}
\newcommand{\erem}{\end{rem}}
\newcommand{\brems}{\begin{rems}}
\newcommand{\erems}{\end{rems}}

\newcommand{\beq}{\begin{eqnarray}}
\newcommand{\beqq}{\begin{eqnarray*}}
\newcommand{\eeq}{\end{eqnarray}}
\newcommand{\eeqq}{\end{eqnarray*}}
\newcommand{\bqn}{\begin{ques}}
\newcommand{\eqn}{\end{ques}}

\newcommand{\vp}{\varphi}

\newcommand{\R}{{\mathbb R}}
\newcommand{\C}{{\mathbb C}}

\newcommand{\D}{{\mathbb D}}

\newcommand{\N}{{\mathbb N}}

\newcounter{minutes}
\divide\time by 60
\newcounter{hours}
\multiply\time by 60 \addtocounter{minutes}{-\time}
\begin{document}
\title[Automorphisms of subalgebras of bounded analytic functions]
{ Automorphisms of subalgebras of bounded analytic functions}

\author{Kanha Behera}
\address{Kanha Behera, Department of Mathematics and Statistics, Indian Institute of Technology, Kanpur - 208016, India.}
\email{beherakanha2560@gmail.com}

\author{Rahul Maurya}
\address{Rahul Maurya, Department of Mathematics and Statistics, Indian Institute of Technology, Kanpur - 208016, India.}
\email{rahulmaurya7892@gmail.com}

\author{P. Muthukumar}
\address{P. Muthukumar,
Department of Mathematics and Statistics, Indian Institute of Technology, Kanpur - 208016, India.}
\email{pmuthumaths@gmail.com, muthu@iitk.ac.in}
%
%


\begin{abstract}
Let $H^\infty$ denote the algebra of all bounded analytic functions on the unit disk. It is well-known that every (algebra) automorphism of $H^\infty$ is a composition operator induced by disc automorphism. Maurya et al., (J. Math. Anal. Appl. 530 : Paper No: 127698, 2024) proved that every  automorphism of the subalgebras $\{f\in H^\infty : f(0) = 0\}$ or $\{f\in H^\infty : f'(0) = 0\}$ is a composition operator induced by a rotation. In this article, we give very simple proof of their results. As an interesting generalization, for any $\psi\in H^\infty$, we show that every automorphism of $\psi H^\infty$ must be a composition operator and characterize all such composition operators. Using  this characterization, we find all automorphism of $\psi H^\infty$ for few choices of $\psi$ with various nature depending on its zeros.
\end{abstract}
\thanks{
File:~\jobname .tex,
          printed: \number\day-\number\month-\number\year,
          \thehours.\ifnum\theminutes<10{0}\fi\theminutes
}
\maketitle
\pagestyle{myheadings}
\markboth{KANHA BEHERA, RAHUL MAURYA and P. MUTHUKUMAR}{Automorphisms of some subalgebras of bounded analytic functions}

\section{Introduction and Preliminaries}\label{Introprelim}
 Let $\D=\{z: |z|<1\}$ be the open unit disc in the complex plane $\C$ and let $\mathbb{T}$ be the (unit circle) boundary of $\D$.
    $\mathrm{Aut}(\D)$ denotes the set of all bijective analytic self-maps (disc automorphisms) of $\D$. For $a\in \D$, $\tau_a$ denotes the special disc automorphism given by $\tau_{a}(z) = (a - z)/(1 - \bar{a}z)$, for $z\in\D$. Define   $\mathcal{H}(\D)$ as the algebra of all analytic functions on $\D$ with point-wise operations.

 A composition operator $C_{\vp}$, associated with an analytic self-map $\vp$ of $\D$, is defined by the following expression:
$$
 	C_{\vp}(f) = f \circ \vp, \ \ \text{for all} \ f \in \mathcal{H}(\D).
$$
It is straightforward to verify that the map $C_{\vp}$ acts as a linear and multiplicative transformation on $\mathcal{H}(\D)$. For non-constant analytic self-map $\vp$ of $\D$, it is trivial to see that $C_\vp$ is injective on $\mathcal{H}(\D)$ using the open mapping theorem and the identity theorem.

Let $H^\infty$ be the Banach algebra of all bounded analytic functions on $\D$, equipped with the supremum norm defined by
$$
\|f\|_\infty = \sup\{|f(z)|: z\in\D\},  \ \ \ \  \ f \in H^\infty.
$$

%

%
 For a given algebra $\mathcal{A}$, by an algebra automorphism of $\mathcal{A}$, we mean a bijective, linear and multiplicative
 self-map of $\mathcal{A}$. In this paper, automorphism of subalgebra of $H^\infty$ will always refer to an algebra automorphism.  For $\vp\in \mathrm{Aut}(\D)$, $C_\vp$ is always an automorphism of $H^\infty$. Conversely, if $T$ is an  automorphism of $H^{\infty}$, then by \cite[Lemma 4.2.1]{Monograph} there is an $\vp\in\mathrm{Aut}(\D)$ such that $T = C_\vp$ i.e.,
 $$
 \text{Automorphisms of} \ H^\infty = \{C_\vp : \vp\in\mathrm{Aut}(\D)\}.
$$
 The result mentioned above can be viewed as a special case of a more general result \cite[Theorem 9]{Rudin},
  by considering  the maximal domain as $\D$.
  Let $A(\D)$ denotes the disk algebra consisting of all continuous functions on $\overline{\D}$ which are analytic in $\D$. With a similar proof as in \cite[Lemma 4.2.1]{Monograph}  we deduce the following result.
\bprop\label{autoAD}
$$
 \text{Automorphisms of} \ A(\D) = \{C_\vp : \vp\in\mathrm{Aut}(\D)\}.
$$
\eprop

  Recently, Maurya et al. proved that every automorphism \cite[Theorem 2.3]{Jaydeb} of $H_0^\infty = \{f\in H^\infty: f(0)=0\}$ or automorphism
   of $H_1^\infty = \{f\in H^\infty: f'(0)=0\}$ \cite[Theorem 3.1]{Jaydeb} is a composition operator induced by a rotation. This naturally leads us to ask the following question:

\textit{Are all automorphisms of any subalgebra  of $H^\infty$ necessarily composition operators induced by some disc automorphisms?}

 As a main result of this article, we prove that
\textit{for any $\psi\in H^\infty(\text{or }A(\D))$, $T$ is an automorphism of $\psi H^\infty(\text{or }\psi A(\D))$ if and only if $T = C_\vp$ for some $\vp\in \mathrm{Aut}(\D)$ such that $\psi \circ \vp = \psi g$, where $g$ is an invertible element of $H^\infty(\text{or }A(\D))$} (See Theorem \ref{thm1} and Corollary \ref{thm1diskalg}).

In this article, we consider general subalgebras of $H^\infty$, not restrict to only closed subalgebras. It is interesting to note that the subalgebra $\vp H^\infty$ of $H^\infty$ is closed if and only if $|\vp|$ is essentially bounded away from zero on the unit circle (see
\cite[Proposition 3.2]{Closed} and remarks below).

This article is organized as follows. In section \ref{sectionj}, we give simpler proof of all the results in Section 2 
 and Theorem $3.1$ of \cite{Jaydeb}. In section \ref{sectionBH}, we give a characterization of composition operator to be an automorphism of $BH^\infty$, where $B$ is a Blaschke product. In section \ref{chara auto}, for a given $\psi\in H^\infty$, we show that every automorphism of $\psi H^\infty$ is a composition operator induced by certain disc automorphism. In section \ref{appli chara auto}, we give characterization of automorphisms of $\psi H^\infty$,
  where $\psi$ is a polynomial that vanishes only at finitely many points on the closed unit disc. Finally in section \ref{section E}, along with various examples, we give another interesting characterization for $C_\vp$ to be an automorphism of $\psi H^\infty$ in terms of Denjoy-Wolff point when $\psi$ is an atomic singular inner function.

   \section{Simpler Proofs of Results in \cite{Jaydeb}}\label{sectionj}

    In this section, we present an alternative and simpler proof of the results established in sections 2 and 3 of \cite{Jaydeb} and extend these results.  For $a \in \D$, define the subalgebras $\mathcal{A}_a = \{f \in H^{\infty} : f(a) = 0\}$ and $\mathcal{B}_a = \{f \in H^{\infty} : f'(a) = 0\}$. 
     We denote $\mathcal{A}_0$ by $\mathcal{A}$ and $\mathcal{B}_0$ by $\mathcal{B}$ for simplicity. Observe that any function $f \in H^\infty$ can be uniquely expressed as $f = f(0) + (f - f(0))$. Therefore, we can decompose $H^\infty$ as the direct sum of $\C$ and $\mathcal{A}$.  That is, $H^\infty = \C \oplus \mathcal{A}$. Similarly, any function $f\in H^\infty$ can be uniquely expressed as $f = z f'(0) + (f - zf'(0))$. Therefore,  we can decompose $H^\infty$ as the direct sum of $z\C$ and $\mathcal{B}$. That is, $H^\infty = z\C \oplus \mathcal{B}$.

     A function $\psi\in H^{\infty}$ is called an inner function if  the radial limit of $\psi$ satisfies  $|\psi(e^{i\theta})| = 1$ a.e. on $\mathbb{T}$.  For example, every disc automorphism  is an inner function.
    Proof of the following result was lengthy and non-trivial (see \cite[Theorem 2.1]{Jaydeb}). Here we give a much simpler and trivial proof of it.

   \bprop\label{inner1}
   Automorphisms of $H^{\infty}$ preserve inner functions.
\eprop
\bpf
   Let $T$ be an automorphism of $H^{\infty}$. By \cite[Lemma 4.2.1]{Monograph}, there exists an
   $\vp\in \text{Aut}(\mathbb{D})$ such that $T = C_{\vp}$.
    For an inner function $\psi$, clearly $T\psi = \psi\circ\vp$ is also an inner function.
\epf

\bthm\label{Thmjay1}
    $T : \mathcal{A} \rightarrow \mathcal{A}$ is an automorphism if and only if there exists $\theta \in \mathbb{R}$ such that
    $$
    T = C_{e^{i\theta}z}.
    $$
\ethm

\bpf
Suppose $T$ is an automorphism of $\mathcal{A}$. Define $S: H^\infty \to H^\infty (= \C \oplus \mathcal{A})$ by
$$
S(\alpha + f) = \alpha + Tf,\,\,\, \alpha \in \C \mbox{ and } f \in \mathcal{A}.
$$
 We claim that $S$ is an automorphism of $H^\infty$. Clearly, $S|_{\mathcal{A}} = T$.

It is straightforward to verify that $S$ is linear, as $T$ is linear. If $\alpha, \beta \in \C$ and $f, g \in \mathcal{A}$, then
$$
S((\alpha + f)(\beta + g)) = \alpha\beta + T(\alpha g + \beta f + fg)
       = (\alpha + Tf)(\beta + Tg).
$$
This shows that $S$ is multiplicative. Suppose $S(\alpha + f) = 0$, that is, $\alpha + Tf = 0.$ 
 By evaluating at $z = 0$, we get $\alpha = 0$
and thus $Tf = 0$. Consequently, we have $f = 0$, implying that $S$ is injective.

Now, for any $g \in H^\infty$, we know $g - g(0) \in \mathcal{A}$. Since $T$ is an automorphism, there exists an $f \in \mathcal{A}$ such that $Tf = g - g(0)$. It yields that $g = g(0) + Tf = S(g(0) + f)$, which proves $S$ is surjective.

Thus, $S$ is an automorphism of $H^\infty$. By \cite[Lemma 4.2.1]{Monograph}, there exists $\varphi \in \text{Aut}(\mathbb{D})$ such that $S = C_{\vp}$.
In particular, we have $z(=id) \in \mathcal{A}$, and $S(z) = T(z) = \varphi \in \mathcal{A}$.
 As $\varphi \in \text{Aut}(\mathbb{D})$ with $\varphi(0)=0$,  $\varphi(z) = e^{i\theta}z$ for some $\theta \in \mathbb{R}$. Consequently, we conclude that $T = C_{e^{i\theta}z}$.

The converse is straightforward to verify, and this completes the proof of the theorem.
\epf

\bthm\label{Thmjay2}
    $T : \mathcal{B}\rightarrow \mathcal{B}$ an automorphism if and only if there exists $\theta \in \mathbb{R}$ such that
     $$
    T = C_{e^{i\theta}z}.
    $$
\ethm

\bpf
Suppose $T$ is an automorphism of $\mathcal{B}$. Notice that $(T(z^3))^2 = (T(z^2))^3$. If $z_0$ is a zero of $T(z^2)$ with multiplicity $n$, then it is also a zero of $T(z^3)$ with multiplicity $\frac{3n}{2}$. Therefore, we define
$$
\tau = \frac{T(z^3)}{T(z^2)},
$$
which is well-defined, and $\tau \in \mathcal{H}(\mathbb{D})$. Since $(T(z^2))^3 = (T(z^3))^2 = \tau^2(T(z^2))^2$, by the identity theorem, we conclude that $T(z^2) = \tau^2$, $T(z^3) = \tau^3$, and $\tau \in H^\infty$.
 Define the map $S: H^\infty (=z\C \oplus \mathcal{B}) \to H^\infty$ by
 $$S(\alpha z + f) = \alpha\tau + Tf,\,\,\,
 \alpha \in \C \mbox{ and } f \in \mathcal{B}.$$
  Clearly $S|_\mathcal{B} = T$. We claim that $S$ is an automorphism  of $H^\infty$.

 Since $T$ is linear, it is trivial to see that $S$ is linear. We now claim that $\tau \notin \mathcal{B}$, i.e., $\tau'(0) \neq 0$. Suppose, for contradiction, that $\tau \in \mathcal{B}$. Then, there exists some $g \in \mathcal{B}$ such that $Tg = \tau$. This would imply $T(g^2) = \tau^2 = T(z^2)$. Consequently, we obtain $g^2 = z^2$. By the identity theorem, this gives $g = \pm z$. However, in both cases $g \notin \mathcal{B}$, leading to a contradiction.

 Next, assume that $S(\alpha z + f) = \alpha\tau + Tf = 0$. Differentiating, we get $\alpha\tau'(0) + (Tf)'(0) = 0$. Since $\tau'(0) \neq 0$ and $(Tf)'(0) = 0$, it follows that $\alpha = 0$. Therefore, we are left with $Tf = 0$, which implies $f = 0$. Hence, $S$ is injective.
  For $g \in H^\infty$, we have $g - \left( \frac{g'(0)}{\tau'(0)} \right)\tau \in \mathcal{B}$. Since $T$ is surjective, $Tf = g - \left( \frac{g'(0)}{\tau'(0)} \right)\tau$ for some $f \in \mathcal{B}$. Therefore $S\left( \frac{g'(0)}{\tau'(0)}z + f \right) = g$, which shows that $S$ is surjective.
  Now, let $\alpha, \beta \in \C$ and $f, g \in \mathcal{B}$. We have
  $$
  (\alpha z + f)(\beta z + g) = (\alpha g(0) + \beta f(0))z + \alpha\beta z^2 + \alpha z(g - g(0)) + \beta z(f - f(0)) + fg.
  $$
 For $h \in \mathcal{B}$ we have $\tau^2 T(zh) = T(z^3 h) = \tau^3 T(h)$, which implies $T(zh) = \tau T(h)$ by identity theorem. Therefore we get $S((\alpha z + f)(\beta z + g)) = S(\alpha z + f)S(\beta z + g)$, which shows that $S$ is multiplicative.

  Hence, $S$ is an automorphism of $H^\infty$. Then there exists  $\varphi \in \text{Aut}(\D)$ such that $S = C_\varphi$. Therefore, $T(f) = C_\varphi(f)$ for all $f \in \mathcal{B}$. In particular, since the square of the identity function is in $\mathcal{B}$, $\varphi^2$ is its image under $T$. Thus, we have $(\varphi^2)'(0) = 2\varphi(0) \varphi'(0) = 0$. Since $\varphi'(0) \neq 0$, this implies that $\varphi(0) = 0$, and therefore $\varphi$ is a rotation, concluding the proof. The converse is trivial to verify, and the theorem is proved.
\epf

Now, we characterize the automorphisms of the subalgebras $\mathcal{A}_a$ and $\mathcal{B}_a$.

\bprop \label{coroHa}
   For $a\in \D$. Let $X$ denote either $\mathcal{A}_a$ or $\mathcal{B}_a$. Then, $T : X\rightarrow X$ is an automorphism if and only if there exists $\theta \in \mathbb{R}$ such that
    $$
    T = C_\vp, \mbox{ where } \vp=\tau_a \circ e^{i\theta}z\circ\tau_a.
    $$
\eprop

   \bpf
   At first, consider the case $X = \mathcal{A}_a$. Let $T$ be an automorphism of $\mathcal{A}_a$. Define a map $S : \mathcal{A}_a \to \mathcal{A}$ by
   $   Sf = f\circ\tau_a   $.
  It is straightforward to verify that $S$ is linear, multiplicative and bijective.
    It follows that $U=STS^{-1}$ is an  automorphism of $\mathcal{A}$. By Theorem \ref{Thmjay1}, $Uf = f\circ e^{i\theta}z$ for some $\theta\in \mathbb{R}$.
    Thus $    T=S^{-1}US= C_{\tau_a \circ e^{i\theta}z \circ \tau_a}.     $

   Now consider the case $X = \mathcal{B}_a$. Let $T$ be an automorphism of $\mathcal{B}_a$. Define the maps $S : \mathcal{B}_a \to \mathcal{B}$ and $U : \mathcal{B} \to \mathcal{B}$ as above. Clearly, $U$ is an  automorphism of $\mathcal{B}$. By using Theorem \ref{Thmjay2}, the desired result follows.
   \epf

For $a\in\D$ and $n\in\N$, define subalgebras: $$\mathcal{A}_a^n=\{f\in H^\infty : f(a) = f'(a) = \cdots = f^{(n-1)}(a) =0\} \mbox{ and }$$
  $$\mathcal{B}_a^n=\{f\in H^\infty :  f'(a) = f''(a)=
   \cdots = f^{(n)}(a) =0\}.$$
  As $H^\infty = \C\oplus z\C\oplus \cdots \oplus z^{n-1}\C\oplus \mathcal{A}_0^n = z\C\oplus \cdots \oplus z^{n}\C\oplus \mathcal{B}_0^n$,
  with a proof similar to those in Theorem \ref{Thmjay1} and  \ref{Thmjay2}, we characterize all automorphisms of $\mathcal{A}_0^n$ and $\mathcal{B}_0^n$ in the following theorem.

\bthm\label{Thmjaygen}
Fix $n\in\N$. Let $X = \mathcal{A}_0^n ~\text{or}~  \mathcal{B}_0^n$ and let $T : X\rightarrow X$ be a map. Then $T$ is an automorphism of $X$ if and only if there exists $\theta \in \mathbb{R}$ such that
     $$
    T = C_{e^{i\theta}z}.
    $$
\ethm


We have the natural isomorphism $C_{\tau_a}$ from $\mathcal{A}_0^n$ and $\mathcal{B}_0^n$ to $\mathcal{A}_a^n$ and $\mathcal{B}_a^n$, respectively. Now using Theorem \ref{Thmjaygen}, we give the following generalisation of Proposition \ref{coroHa}.

\bprop \label{coroH_an}
    Let $X$ denote either $\mathcal{A}_a^n$ or $\mathcal{B}_a^n$. Then, $T : X\rightarrow X$ is an automorphism if and only if there exists $\theta \in \mathbb{R}$ such that
    $$
    T = C_\vp, \mbox{ where } \vp=\tau_a \circ e^{i\theta}z\circ\tau_a.
    $$
\eprop

For $a\in\D$ and $n\in\N$, define the subalgebras of $A(\D)$: $$\mathcal{\tilde{A}}_a^n=\{f\in A(\D) : f(a) = f'(a) = \cdots = f^{(n-1)}(a) =0\} \mbox{ and }$$
  $$\mathcal{\tilde{B}}_a^n=\{f\in A(\D) :  f'(a) = f''(a)=
   \cdots = f^{(n)}(a) =0\}.$$

  \bcor \label{coroAD_an}
    Let $X$ denote either $\mathcal{\tilde{A}}_a^n$ or $\mathcal{\tilde{B}}_a^n$. Then, $T : X\rightarrow X$ is an automorphism if and only if there exists $\theta \in \mathbb{R}$ such that
    $$
    T = C_\vp, \mbox{ where } \vp=\tau_a \circ e^{i\theta}z\circ\tau_a.
    $$
\ecor
The proof is similar to that of Proposition \ref{coroH_an}.

%
%

   \beg
   Consider the following subalgebras $\mathcal{A} = \{f\in H^\infty:f(0)=0\}$ and $\mathcal{C}=\{f\in H^\infty:f(0)=0=f(1/2)\}$. Then it is easy to see that the composition operator $C_\vp$ induced by the disc automorphism $\tau_{1/2}$ is an algebra automorphism of $\mathcal{C}$ with own inverse i.e., $C_\vp\circ C_\vp = I$. It is trivial to see that $C_\vp$ cannot be extended to (a composition operator induced by rotation) an automorphism  of $\mathcal{A}$.

   \eeg
   Therefore, in general an automorphism of subalgebra may not be extended as an automorphism of the given algebra. However in this section, we have proved that every automorphism of $\mathcal{A}_a^n$ or $\mathcal{B}_a^n$  can be extended as an automorphism of $H^\infty$. Consequently, we conclude that automorphisms of these subalgebras are also composition operators.

    \section{$C_\vp$ as an automorphism of $BH^\infty$}\label{sectionBH}

   In section \ref{sectionj}, we observe that every automorphism of $\mathcal{A}_a^n (=\tau_a^n H^\infty)$ 
    is a composition operator induced by some disc automorphism. Motivated by this result, in this section, we discuss  which composition operators are automorphism of $BH^\infty$, where $B$ is a Blaschke product.

    For a  finite or infinite sequence $\{z_j\}$ in $\D$ with $\sum\limits_{j}(1-|z_j|) < \infty, m\in\N\cup\{0\}$ and $\gamma\in\mathbb{T}$,
    the product
   $$
   B(z) = \gamma z^m\prod_{j}\frac{|z_j|}{z_j}\frac{z_j - z}{1 - \bar{z_j}z}, \ \ \ \ \ z\in\D
   $$
   is called Blaschke product. It is well-known that $B$ is an analytic function with $|B(z)|<1$ on $\D$ and $|B(e^{i\theta})|= 1$ a.e.
    (see \cite[Theorem 2.4]{Duren}). As any change in the unimodular constant $\gamma$ does not alter the algebra $BH^\infty$, without loss
    of generality, we may assume that $\gamma=1$.

   We give the statement of the following result (\cite[Theorem 2.5]{Duren}) of F. Riesz, which we have used in the proof of Lemma \ref{multzero1}, which, in turn, is used throughout the section:
Every function $f(z) \not\equiv 0$ from the class $H^p$ ($0<p \leq\infty$) can be factored in the form $f(z) = B(z)g(z)$, where $B(z)$ is a Blaschke product and $g(z)$ is an $H^p$ function which does not vanish in $|z|<1$.

For $f\in\mathcal{H}(\D)$, $Z(f)$ denotes the set of all zeros of $f$ inside $\D$. For $w\in Z(f)$, we denote its multiplicity by $mult_f(w)$.

    \blem\label{multzero1}
    Let $\vp \in \text{Aut}(\D)$ and $\psi\in H^\infty$. Then  $mult_{\psi\circ \vp}(w) = mult_\psi(\vp(w))$ for all $w\in Z(\psi\circ\varphi)$.
    \elem

    \bpf
    Suppose $w\in Z(\psi \circ \varphi)$ and $mult_{\psi \circ \varphi}(w) = k$. As a consequence of \cite[Theorem 2.5]{Duren}, we have
    $
    \psi \circ \varphi = \tau_w^{k} g_1,
    $
    where $g_1 \in H^\infty$ and $g_1(w) \neq 0$. Since $\varphi$ is a disc automorphism, the above relation implies
    $
    \psi = (\tau_w \circ \varphi^{-1})^{k} \left( g_1 \circ \varphi^{-1} \right).
    $
    Also, we have $\tau_w \circ \varphi^{-1} = \alpha \tau_{\varphi(w)}$ for some $\alpha \in \mathbb{T}$, because $\varphi(w)$ is a zero of disc automorphism $\tau_w \circ \varphi^{-1}$.
    Since $\varphi(w)$ is not a zero of $g_1 \circ \varphi^{-1}$, we conclude that $mult_\psi(\varphi(w)) = k$, which implies that
    $$
    mult_{\psi \circ \varphi}(w) = mult_\psi(\varphi(w)).
    $$
    \epf

    The following remark shows that the equality $mult_{\psi \circ \vp}(w) = mult_\psi(\vp(w))$ does not necessarily hold for a general $\varphi$.

    \brem
    Consider the functions $\psi(z) = z$ and $\varphi(z) = z^2$. In this case, we have $mult_{\psi \circ \vp}(0) \neq mult_\psi(\vp(0))$.
    \erem
    For a  self-map $\vp$ of $\D$, we denote $\vp \circ\vp \circ \cdots \circ \vp$  $(n \ \text{times})$  by $\vp_n$.

     \bthm\label{finiteblascke}
     Let $\vp$ be an analytic self-map on $\D$ and $B$ be a finite Blaschke product. Then $C_\vp : BH^\infty \to BH^\infty$ is an automorphism if and only if $\vp \in \mathrm{Aut}(\D)$ and $mult_{B}(\vp(w)) = mult_B(w)$ for all $w\in Z(B)$.
     \ethm
     \bpf
     If $Z(B)$ is a singleton, then the claimed result follows directly from the conclusions of the previous section.

     Now, let $n \geq 2$ and  $Z(B) = \{w_1, w_2, \ldots, w_n\}$. Thus, $B = \tau_{w_1}^{k_1} \tau_{w_2}^{k_2} \cdots \tau_{w_n}^{k_n}$, with each $k_i$ is a natural number. Suppose $C_\varphi$ is an automorphism of $BH^\infty$. Then $\vp$ is non-constant and $B \circ \varphi = Bg$ for some $g \in H^\infty$.
    Therefore  $\varphi(Z(B)) \subset Z(B)$. We claim that $\varphi$ maps $Z(B)$ to itself bijectively.
    Suppose, for contradiction, that $\varphi$ is not surjective. Then, there exists some $w_j \in Z(B) \setminus \varphi(Z(B))$. Define a function $g$ by
     $$
     g= \prod_{i=1,  i \neq j}^{n} \left( \tau_{w_i} \circ \varphi \right)^k= \prod_{i=1,  i \neq j}^{n} \left( \tau_{w_i} \right)^k \circ \varphi,
     $$
     where $k = \max\{k_1, k_2, \dots, k_n\}$. It is straightforward to see that 
     $g \in BH^\infty$. Since $C_\varphi$ is an automorphism, there exists $f \in BH^\infty$ such that $C_\varphi f = g$.
    %
    Since $C_\varphi$ is injective on $H^\infty$, we have
    $$
     f = \prod_{i=1, i \neq j}^{n} (\tau_{w_i})^k.
     $$
     This leads to $f(w_j) \neq 0$, implying that $f \notin BH^\infty$, which is a contradiction. Therefore, $\varphi(Z(B)) = Z(B)$. Since $Z(B)$ is finite, $\varphi$ maps $Z(B)$ to itself bijectively. 
     Hence, $\varphi_n$ fixes $w_i$ for all $i = 1, 2, \dots, n$. Since the self-map $\varphi_n$ of $\D$ has more than one fixed point in $\D$, $\varphi_n$ must be  identity. Consequently, $\varphi\in\text{Aut}(\mathbb{D})$.

      Now we claim that $mult_{B}(\varphi(w)) = mult_B(w)$ for all $w \in Z(B)$. Let $w_s \in Z(B)$. If $\varphi(w_s) = w_s$, then our claim is immediately true. Suppose instead that $\varphi(w_s) = w_t$ for some $t \neq s$. Suppose $k_s < k_t$.
      Define
     $$
     h = (\tau_{w_t} \circ \varphi)^{k_s} \prod_{i=1, i \neq t}^{n} (\tau_{w_i} \circ \varphi)^k
     $$
     where $k = \max \{k_1, k_2, \ldots, k_n\}$. It is clear that $h \in BH^\infty$. But, 
     $$
     C_{\varphi}^{-1}(h) = \tau_{w_t}^{k_s} \prod_{i=1, i \neq t}^{n} \tau_{w_i}^k
     $$
     which does not belong to $BH^\infty$, which is a contradiction. 
      Similarly, we derive a contradiction if $k_s > k_t$. Thus, we must have $mult_{B}(\varphi(w)) = mult_B(w)$ for all $w \in Z(B)$.

      For the converse part, suppose that $mult_{B}(\vp(w)) = mult_B(w)$ for all $w\in Z(B)$ and $\vp\in \text{Aut}(\D)$. By Lemma \ref{multzero1}, $mult_{B\circ\vp}(w) = mult_B(w)$. So, $B\circ\vp =\alpha B$ for some $\alpha\in\mathbb{T}$, because
       both $B\circ\vp$ and $B$ are  Blaschke products with the same zeros of the same multiplicities. It implies that $BH^\infty$ is invariant under both $C_\varphi$ and $C_{\vp^{-1}}$. Thus, $C_\varphi$ is an automorphism on $BH^\infty$, which completes the proof of the theorem.
      \epf

     \bcor \label{anyblaschke}
     Let $\vp \in\mathrm{Aut}(\D)$ and $B$ be an Blaschke product. Then $C_\vp$ is an automorphism of $BH^\infty$ if and only if $mult_{B}(\vp(w)) = mult_B(w)$ for all $w\in Z(B)$.
     \ecor
     \bpf Assume that $C_\vp$ is an automorphism of $BH^\infty$. It gives that $B \circ \varphi = Bg$ for some $g \in H^\infty$. So, by Lemma \ref{multzero1}
     $mult_B(w)\leq mult_{B g}(w)= mult_B(\vp(w))$ for all $w\in Z(B).$
     As $C_{\varphi^{-1}}$ is also an automorphism of $BH^\infty$, similarly we have,
     $$
     mult_B(w) \leq mult_B(\varphi^{-1}(w)) \quad \text{for all} \quad w \in Z(B).
     $$
    Since $\vp(w)\in Z(B)$ for all $w\in Z(B)$, we get
     $$
     mult_B(\varphi(w)) = mult_B(w)  \quad \text{for all} \quad w \in Z(B).
     $$

     The converse part can be verified in the same way as in Theorem \ref{finiteblascke}. Here is an alternate proof of it.
      Suppose $mult_{B}(\vp(w)) = mult_B(w)$ for all $w\in Z(B)$. By Lemma \ref{multzero1},  we get    $mult_{B\circ\vp}(w) = mult_B(w)$. As noted above Corollary 1.4 in \cite{Pmuthu}, $BH^\infty$ is invariant under both $C_\varphi$ and $C_{\varphi^{-1}}$. Thus, $C_\varphi(BH^\infty)=BH^\infty$ and therefore $C_\varphi$ is an automorphism on $BH^\infty$.

           \epf

     \brem
     Corollary \ref{anyblaschke} need not hold for a general subalgebra $\psi H^\infty$, where $\psi \in H^\infty$. For example,
       let $\psi= 1 + z$ and $\vp= \tau_{1/2}$. Since $Z(\psi)$ is empty, $mult_{\psi}(\vp(w)) = mult_{\psi}(w)$ for all $w\in Z(\psi)$ holds trivially. As $C_{\vp}\psi\notin \psi H^\infty$, $C_\vp$ cannot be an automorphism of $\psi H^\infty$.
     \erem

   \section{Automorphisms of $\psi H^\infty$}\label{chara auto}

   In this section, for any $\psi\in H^\infty$ we prove that every automorphism of $\psi H^\infty$
    can be extended to an automorphism of $H^\infty$ and hence it is a composition operator. Moreover, we characterize all such composition operators.  Throughout the article, we assume, without loss of generality, that $\psi \in H^\infty$ is a non-identically zero function, i.e., $\psi \not\equiv 0$ on $\mathbb{D}$.
   Recall that a function $f\in H^\infty$ is said to be invertible if there exists  $g\in H^\infty$ such that $fg \equiv 1$ on $\D$, that is $1/g \in H^\infty$.

\blem \label{lem1}
      Let $T$  be an automorphism of $\psi H^\infty$. Then $T(\psi) = \psi g$, for some invertible function $g\in H^\infty$.
\elem

\bpf
     Since $T$ is an automorphism of $\psi H^\infty$ and $T\psi \in \psi H^\infty$, there exists some $g\in H^{\infty}$ such that $T\psi = \psi g$. Similarly, $T^{-1}(\psi) = \psi g_1$ for some  $g_1 \in H^\infty$. By the multiplicativity of $T^{-1}$, we have  $T^{-1}(\psi^2) = (T^{-1}(\psi))^2 = \psi^2g_1^2$. This implies that
     $$
       \psi^2 = T(\psi^2g_1^2) = T(\psi)T(\psi g_1^2) = (\psi g)(\psi g_2),  \ \ \ \text{for some} \ g_2\in H^\infty.
     $$
     Therefore, we have $\psi^2(gg_2 - 1) = 0$. As $\psi\not\equiv
      0$, by the identity theorem, we get $gg_2 \equiv 1 $. Thus $g$ is an invertible in $H^\infty$. It completes the proof.
\epf
Note that the converse of Lemma \ref{lem1} is not true. That is,  even if $T\psi=\psi g$ for some invertible function $g\in H^\infty$, $T$ may not be an algebra automorphism of $\psi H^\infty$. For example, consider $\psi\in H^\infty$ to be arbitrary and $Tf = 2f$ for all $f\in \psi H^\infty$. Clearly, $T$ is not multiplicative on $\psi H^\infty$.

    \bprop \label{pro1}
   Suppose $\psi\in H^\infty$ and $\vp\in \text{Aut}(\mathbb{D})$ such that $\psi\circ\vp = \psi g$ for some invertible function $g\in H^\infty$. Then the composition operator $C_\varphi$ is an automorphism of $\psi H^\infty$.
   \eprop

   \bpf
   Since $\psi\circ\vp = \psi g $, it is trivial to see that $C_\vp(\psi H^\infty) \subseteq \psi H^\infty$. As $\vp$ is non-constant, $C_\varphi$ is injective. To see $C_\varphi$ is onto, fix $\psi g_1 \in \psi H^\infty$. Consider $f = (g_1/g) \circ \vp^{-1}$, which is clearly in $H^\infty$. Then $C_{\vp}(\psi f) = \psi g_1$. As $C_\varphi$ is always linear and multiplicative, we get $C_\varphi$ is an automorphism of $\psi H^\infty$.
   \epf

\bthm \label{thm1}
    Fix $\psi\in H^\infty$. Then $T$ is an automorphism of $\psi H^\infty$ if and only if $T = C_\vp$ for some $\vp\in \mathrm{Aut}(\D)$ such that $\psi \circ \vp = \psi g$, where $g$ is an invertible element of $H^\infty$.
\ethm

\bpf
    Let $T$ be an automorphism of $\psi H^\infty$. Define $S : H^\infty \to H^\infty$ by
     $$
     Sf = \frac{T(\psi f)}{T\psi}, \,\, f\in H^\infty.
     $$
      By  Lemma \ref{lem1}, the map $S$ is well-defined and  agrees with $T$ on $\psi H^\infty$. Now,
       we claim that $S$ is an automorphism $ H^\infty$. As $T$ is linear, $S$ is linear. If $f, g\in H^\infty$ then
     $$
     S(fg)  = \frac{T(\psi fg)}{T\psi} = \frac{T(\psi^2fg)}{T\psi\cdot T\psi} = \frac{T(\psi f)}{T\psi}\cdot\frac{T(\psi g)}{T\psi} = Sf \cdot Sg,
     $$
     which shows that $S$ is multiplicative.

     Suppose $Sf = 0$ for some $f\in H^\infty$. This implies that $T(\psi f) = 0$. Since $T$ is an automorphism $\psi f = 0$. As $\psi \not \equiv 0 $, identity theorem forces that $f \equiv 0$. Thus $S$ is injective. Consider $g\in H^\infty$. Since $gT\psi \in \psi H^\infty$ and $T^{-1}$ is an automorphism of $\psi H^\infty$, we get $T^{-1}(gT\psi) \in \psi H^\infty$. Therefore, $S\left( \frac{T^{-1}(gT\psi)}{\psi} \right) = g$, which proves $S$ is onto. Thus $S$ is an automorphism of $H^\infty$. By \cite[Lemma 4.2.1]{Monograph}, there exists $\vp \in \text{Aut}(\D)$ such that $S = C_\vp$. Hence $Tf = C_\vp f$, for all $f \in \psi H^\infty$. By Lemma \ref{lem1}, $T\psi = \psi \circ \vp = \psi g$ for some invertible element $g\in H^\infty$. Converse follows from Proposition \ref{pro1}. This completes the proof of the theorem.
\epf

 \brem
 In general, automorphisms of a Banach algebra need not be an isometry. For any $\vp\in\text{Aut}(\D)$, we have  $\|C_\vp f\|_{\infty} = \|f\circ\vp\|_\infty = \|f\|_\infty$ for all $f\in \psi H^\infty$, hence every automorphism of $\psi H^\infty$ is  always an isometry.
 \erem

 With a proof almost identical to that of Theorem \ref{thm1}, we can derive the following result.
 \bcor\label{thm1diskalg}
 Fix $\psi\in A(\D)$. Then $T$ is an automorphism of $\psi A(\D)$ if and only if $T = C_\vp$ for some $\vp\in \mathrm{Aut}(\D)$ such that $\psi \circ \vp = \psi g$, where $g$ is an invertible element of $A(\D)$.
 \ecor
 \brem\label{rmkpsiAD}
 Whenever $\psi\in A(\D)$, results on $\psi H^\infty$ in this article continue to hold on $\psi A(\D)$. Thus, we do not make separate statements for the latter case of $\psi\in A(\D)$.
 \erem

\bcor \label{inner2}
    Automorphisms of $\psi H^{\infty}$ preserve inner functions.
\ecor

\bcor \label{inner3}
     Let $\psi$ be an inner function and $T$ be an automorphism of $\psi H^\infty$. Then there exists an $\alpha\in \mathbb{T}$ such that $T(\psi) = \alpha\psi$.
\ecor

\bpf
  By Lemma \ref{lem1}, $T\psi = \psi g$ for some invertible element $g\in H^\infty$. Since $\psi$ is inner by Corollary \ref{inner2}, $\psi g$ is also an inner function leading to the equality $|g| =|\psi g| = 1$ a.e on $\mathbb{T}$. Thus $g$ and $1/g$ are inner functions. This yields that $|g(z)|= 1$ on $\D$. Consequently, $g$ is a unimodular constant.
\epf

   Using our main theorem we now give an alternative proof of Corollary \ref{anyblaschke}.

  \bcor 
  Let $\vp\in \mathrm{Aut}(\D)$ and  $B$ be Blaschke product. Then $C_\vp$   is an automorphism of $BH^\infty$ if and only if $mult_{B}(\vp(w)) = mult_B(w)$ for all $w\in Z(B)$.
  \ecor
  \bpf
  Let $C_\vp$   be an automorphism of $BH^\infty$.  According to Theorem \ref{thm1}, we have $B \circ \varphi = B g$ for some invertible element $g$ in $H^\infty$. This implies that $mult_{B \circ \varphi}(w) = mult_{Bg}(w) = mult_B(w)$ for all $w \in Z(B)$. Since $mult_{B\circ\varphi}(w) = mult_B(\varphi(w))$ for all $w \in Z(B)$, we conclude that $mult_B(\varphi(w)) = mult_B(w)$ for all $w \in Z(B)$.
   The converse part is already verified in Corollary \ref{anyblaschke}.
  \epf

  The proof of the following result is exactly the same as in the above corollary.
  \bcor \label{multzeropsi}
  Let $\vp\in \mathrm{Aut}(\D)$ and  $\psi\in H^\infty$. If $C_\vp$  is an automorphism of $\psi H^\infty$, then $mult_{\psi}(\vp(w)) = mult_\psi(w)$ for all $w\in Z(\psi)$.
  \ecor

  Recall that $\vp\in \mathrm{Aut}(\D)$ other than identity is called elliptic if $\vp$ has a fixed point in $\D$.

\bprop \label{elliptic}
      Let $T:\psi H^\infty \to \psi H^\infty$ be an automorphism, where $\psi$ is a nonconstant inner function, which extends continuously up to $\overline{\D}$. Then  $T = C_\vp$ for some $\vp$, which is identity or  an elliptic automorphism of $\D$.
\eprop

\bpf
      By Theorem \ref{thm1}, we have  $T = C_\vp$, for some $\vp\in Aut(\D)$. It remains to show that $\vp$ is elliptic. Since $\psi$ is an inner function, Corollary \ref{inner3} yields $T\psi = \alpha\psi$, for some $|\alpha|=1$. This implies, $ T^n\psi = \psi\circ\vp_n = \alpha^n\psi$. Assume, for the sake of contradiction, that $\vp$ is neither elliptic nor identity. Then, by Denjoy-Wolff Theorem \cite[Section 5.1]{Shapiro}, $n$-fold composition $\vp_n$ converges uniformly to a point $p\in \partial{\mathbb{D}}$ on compact subset of $\D$. Consequently,  $\psi(\vp_n(z))$ converges to $\psi(p)$ for each $z\in \D$. Thus, $\lim_{n \to \infty} \alpha^n\psi(z) = \psi(p)$, implying  $|\psi(z)| = |\psi(p)| = c$, for all $z\in \D$ and thus  $\psi$ is constant function.  This contradicts the assumption that $\psi$ is nonconstant. Therefore, $\vp$  must be identity or an elliptic automorphism of $\D$.
\epf

For any $\psi \in H^\infty$, Theorem \ref{thm1} asserts that characterizing the automorphisms of $\psi H^\infty$ is equivalent to characterizing the composition operators $C_\varphi$ that act as automorphisms of $\psi H^\infty$, where $\varphi \in \mathrm{Aut}(\mathbb{D})$. Accordingly, in Sections $5$ and $6$, we will analyze the composition operator $C_\varphi$ for a given $\varphi \in \mathrm{Aut}(\mathbb{D})$, focusing on determining the conditions under which it becomes an automorphism of $\psi H^\infty$.

   \section{Special Case}\label{appli chara auto}

   In this section, we restrict our attention to the special case that $\psi$ has only finitely many zeros in the closed unit disc $\overline{\D}$.

    \bprop\label{psiBg}

   Suppose that $\psi \in H^\infty$ such that $\psi = Bg$, where $B$ is a Blaschke product and $g \in H^\infty$ is non-vanishing on $\mathbb{D}$
   and let $\varphi \in \mathrm{Aut}(\D)$. Then $C_\varphi$ is an automorphism of $\psi H^\infty$ if and only if $C_\varphi$ is an automorphism of both $B H^\infty$ and $g H^\infty$.
   \eprop

   \bpf
   Let $C_\varphi$ be an automorphism of both $BH^\infty$ and $gH^\infty$. Then, by Lemma \ref{lem1}, we have $B \circ \varphi = B g_1$ and $g \circ \varphi = g g_2$, for some invertible elements $g_1$ and $g_2$ of $H^\infty$. Therefore,
   $$
   \psi \circ \varphi = (B \circ \varphi)(g \circ \varphi) = \psi g_1 g_2,
   $$
   where $g_1 g_2$ is an invertible element of $H^\infty$. Hence, by Proposition \ref{pro1}, $C_\varphi$ is also an automorphism of $\psi H^\infty$.

   Conversely, suppose that $C_\vp$ is an automorphism of $\psi H^\infty$.
   Again by Lemma \ref{lem1}, we have $\psi \circ \varphi = \psi g_3$ for some invertible elements $g_3$.
    By Corollary \ref{multzeropsi}, it follows that $mult_\psi(\vp(w)) = mult_\psi(w)$ for all $w \in Z(\psi)$. Since $\psi = Bg$ with $g$ being non-vanishing on $\mathbb{D}$, we deduce that
   $mult_B(w) = mult_B(\vp(w))$ for all $w \in Z(B)$. By Corollary \ref{anyblaschke}, it follows that $C_\vp$ is an automorphism of $B H^\infty$.

   Furthermore, Corollary \ref{inner3} implies that $B \circ \vp = \alpha B$, where $\alpha \in \mathbb{T}$. By Theorem \ref{thm1}, we have
   $$
   Bg g_3 = \psi g_3 = \psi \circ \vp = (B \circ \vp)(g \circ \vp) = \alpha  (g \circ \vp) B.
   $$
 By identity theorem, it follows that $g \circ \varphi = g(\overline{\alpha} g_3)$. Consequently, by Proposition \ref{pro1}, $C_\varphi$ is an automorphism of $g H^\infty$.
   \epf

We will now provide a characterization of the composition operator in terms of the multiplicity of zeros of $\psi$ when it acts as an automorphism of $\psi H^\infty$, where $\psi$ is a polynomial with roots on $\mathbb{T}$. Before proceeding, we will establish a lemma that serves as a crucial tool in proving this result. The following lemma uses the fact that if $\vp\in \text{Aut}(\D)$ then $\vp$ is analytic on a domain containing the closed unit disk. If $\vp\in \text{Aut}(\D)$ then
$\vp(z)=\eta\frac{a-z}{1-\overline{a}z}$ for some $a\in\D$ and $|\eta| =1$. As $1/\overline{a}$, the only singularity of $\vp$, lies outside of $\overline{\D}$, we get $\vp$ is analytic on a domain containing the closed unit disk.

   \blem\label{autominus}
   If $\vp\in\text{Aut}(\D)$ such that $\vp(a) = b$ for $a,b\in\overline{\D}$, then $\vp(z)- b = (z-a)g(z)$ where $g\in H^\infty$ such that $1/g$ is also in $H^\infty$.
   \elem
   \bpf
   Since $\vp$ is analytic on $\overline{\D}$ and $\vp(a) = b$, we have $\vp(z)- b = (z-a)g(z)$ where $g$ is analytic on $\overline{\D}$.
    Then $g(a) = \vp'(a) \neq 0$ and $g$ does not vanish also on $\overline{\D}\setminus\{a\}$ because $\vp$ is injective.
    It forces that  $1/g$ is in $H^\infty$.
   \epf

   \bthm\label{Poly1}
   Let $\vp\in \mathrm{Aut}(\D)$ and let $\psi(z) = (z - w_1)^{n_1} (z - w_2)^{n_2} \cdots (z - w_k)^{n_k}$ be a polynomial,
   where  $|w_j|=1$ and  $n_j\in \mathbb{N}$ for all $j$.
   Then $C_\varphi$ is an automorphism of  $\psi H^\infty$ if and only if
   $$
   mult_{\psi}(\varphi(w_j)) = mult_{\psi}(w_j),  \ 1 \leq j \leq k.
   $$
     \ethm
   \bpf
   Suppose that $mult_{\psi}(\varphi(w_j)) = mult_\psi(w_j)$ for $1 \leq j \leq k$, and let $\varphi \in \mathrm{Aut}(\mathbb{D})$. This implies that $\varphi$ bijectively maps $\{w_i:1 \leq i \leq k\}$ to itself. Therefore, for each
    $j$, there exists $l_j$ such that $\varphi(w_j) = w_{l_j}$.
     Moreover, we can express $\psi$ in the form
   $$
   \psi(z) = (z - w_{l_1})^{n_1} (z - w_{l_2})^{n_2} \cdots (z - w_{l_k})^{n_k}.
   $$

   By Lemma \ref{autominus}, we have
%
   $$\psi\circ\vp(z)  = (\vp(z) - w_{l_1})^{n_1} (\vp(z) - w_{l_2})^{n_2} \cdots (\vp(z) - w_{l_k})^{n_k}=
     \psi(z)g(z), $$

   for some invertible $g$ in $H^\infty$.
   Therefore by Proposition \ref{pro1}, $C_\vp$ is an automorphism of $\psi H^\infty$.

    Conversely, suppose $C_\vp$ is an automorphism of $\psi H^\infty$. Then by Lemma \ref{lem1}, we get $\psi\circ\vp = \psi g$, for some invertible element $g\in H^\infty$. As $\psi\circ\vp$ is analytic on $\overline{\D}$ and  $g\in H^\infty$, there exists $M>0$ with
    $$
    |\psi\circ\vp(z)|\leq M|\psi(z)| \ \ \text{on } \D, \  \ \psi \circ \varphi(w_j) = 0 \ \text{for} \  1 \leq j \leq k.
    $$
  This gives that for each $j$, $\varphi(w_j)$ is a zero of $\psi$ and thus $\varphi$ maps $ \{w_i:1 \leq i \leq k\}$ to itself bijectively, as $\vp$ is 1-1.
  Now we claim that $mult_\psi(\varphi(w_j)) = mult_\psi(w_j)$ for $1 \leq j \leq k$.
   Fix $1 \leq s \leq k$. If $\varphi(w_s) = w_s$, the claim holds trivially. Now, suppose  $\varphi(w_s) = w_t$ for some $t \neq s$ and  assume $n_s > n_t$ for a contradiction. Then  by Lemma \ref{autominus}, there exists an invertible $g_s\in H^\infty$ such that

   \begin{align*}
   \psi \circ \vp(z)  
   & = (\vp(z) - w_t)^{n_t}\prod_{i = 1, i\neq t}^{k} (\vp(z) - w_i)^{n_i} \\
   & = (z - w_s)^{n_t}(g_s(z))^{n_t}\prod_{i = 1, i\neq t}^{k} (\vp(z) - w_i)^{n_i}.
   \end{align*}
  Also we have $\psi\circ\vp(z) =\psi(z)g(z)$ 
  for some invertible  $g$ in $H^\infty$ and therefore
  $$
  (z - w_s)^{n_t}(g_s(z))^{n_t}\prod_{i = 1, i\neq t}^{k} (\vp(z) - w_i)^{n_i}  =
  (z - w_s)^{n_s}g(z)\prod_{i = 1, i\neq s}^{k} (z - w_i)^{n_i}.
  $$
  Since $g$ and $1/g_s$ are in $H^\infty$, there exists $N>0$ such that
 \begin{align}\label{eqbd}
  \prod_{i = 1, i\neq t}^{k} |\vp(z) - w_i|^{n_i} &\leq
     N|z - w_s|^{n_s - n_t}\prod_{i = 1, i\neq s}^{k} |z - w_i|^{n_i} \mbox{ on } \D.
     \end{align}
  Letting $z\rightarrow w_s$ in \eqref{eqbd}, we have
  $$
  \prod_{i = 1, i\neq t}^{k} (\vp(w_s) - w_i)^{n_i} = 0,
  $$
 which is a contradiction as $\varphi(w_s) = w_t$. A similar contradiction occurs if $n_s < n_t$. Therefore, we must have $n_s = n_t$.
  This completes the proof.
   \epf

\bcor\label{zeroinDclos}
Let $\psi = B g \in H^\infty$, where $B$ is a Blaschke product and $g(z) = (z - w_1)^{n_1} (z - w_2)^{n_2} \cdots (z - w_k)^{n_k}$ where  $|w_j|=1$ and  $n_j\in \mathbb{N}$ for all $j$.Then for $\vp\in \mathrm{Aut}(\D)$, $C_\vp$ is an automorphism of $\psi H^\infty$ if and only if $mult_B{\vp(w)} = mult_B{w}$ for all $w\in Z(B)$ and $mult_g{\vp(w_j)} = mult_g{w_j}$ for $1 \leq j \leq k$.
\ecor
\bpf
The proof follows from Proposition \ref{psiBg}, Theorem \ref{Poly1} and Corollary \ref{anyblaschke}.
\epf

   \bcor
  Let $\psi = B g$, where $B$ is a finite Blaschke product and $g$ as in Corollary \ref{zeroinDclos}.
   Assume that the multiplicity of one of the zeros of $g$ is different from the multiplicities of the other zeros. Then, the identity automorphism is the only automorphism of the algebra $\psi H^\infty$.

   \ecor
   \bpf
 Let $C_\vp$ be an automorphism of $\psi H^\infty$. Corollary \ref{psiBg} establishes that $C_\vp$ is also an automorphism of both $BH^\infty$ and $gH^\infty$. Consequently, by Proposition \ref{elliptic}, the mapping $\vp$ must be either an elliptic  or the identity. In addition, Theorem \ref{Poly1} implies that $\vp$ fixes a zero of $g$. This means that $\vp$ cannot be elliptic, so it must be the identity map on $\D$. Therefore, any automorphism of $\psi H^\infty$ is  trivial.
  \epf

   \bcor
   Let  
   $\psi\in H^\infty$ such that there exist $a,b\in Z(\psi)$ with $mult_\psi(a)=m_1\neq m_2=mult_\psi(b)$ and the multiplicity of any other zeros of
   $\psi$ is neither $m_1$ nor $m_2$. Then the only automorphism of $\psi H^\infty$ is identity.
       \ecor

   \bpf
   Let $C_\vp$ be an automorphism of $\psi H^\infty$. 
   Then, by Corollary \ref{multzeropsi}, it follows that $\phi(a) = a$  and $\phi(b) = b$. As $\vp$ has two fixed points in $\D$, $\vp$ must be identity on $\D$. It completes the proof.
   \epf

   \section{Examples}\label{section E}
   In this section, for  various special choices of $\psi$ with different nature depending on its zeros, we provide simpler characterization for automorphisms of the subalgebra $\psi H^\infty$
    and hence determine all automorphisms of $\psi H^\infty$.
    %

      \bthm \label{thm58}
   Let $\vp\in \mathrm{Aut}(\D)$ and $\psi(z) = \exp({\alpha\frac{z + 1}{z - 1}}), \alpha>0$. Then $C_\vp$ is an automorphism of $\psi H^\infty$ if and only if $\vp(1) = 1$ and $\vp^{'}(1) = 1$ i.e., $\vp$ is a parabolic disc automorphism
   with fixed point $1$.
   \ethm
   \bpf
   Suppose $\varphi \in \mathrm{Aut}(\mathbb{D})$ such that $\varphi(1) = 1$ and $\varphi'(1) = 1$. Then by using
   \cite[Theorem 2.44]{Cowen-Book} and Julia Lemma \cite[Lemma 2.41]{Cowen-Book} for $\vp$ and $\vp^{-1}$, we get
   $$
   \frac{|1 - \vp(z)|^2}{1 - |\vp(z)|^2} = \frac{|1 - z|^2}{1 - |z|^2}
   $$
   That is,
   $$
   \text{Re}\bigg(\frac{\vp(z) + 1 }{\vp(z) - 1}\bigg) = \text{Re}\bigg(\frac{z + 1 }{z - 1}\bigg).
   $$
   It leads that
   \begin{equation}\label{e1singular}
   	\frac{\vp(z) + 1}{\vp(z) - 1} - \frac{z + 1}{z - 1} = ic,
   \end{equation}
   for some $c\in \mathbb{R}$. From the equation (\ref{e1singular}), we have
   \begin{align*}
   	\exp\left({\alpha\frac{z + 1}{z - 1}}\right)\circ \vp = \exp\left({\alpha\frac{\vp(z) + 1}{\vp(z) - 1}}\right)
   	  = \exp\left({\alpha\frac{z + 1}{z - 1}}\right) \exp(i\alpha c).
   \end{align*}
    By Proposition \ref{pro1},
    $C_\vp$ is an automorphism of $e^{\alpha\frac{z + 1}{z - 1}}H^\infty$.

   Now, we prove the forward direction. Let $C_\vp$ be an automorphism of $e^{\alpha\frac{z + 1}{z - 1}}H^\infty$. By Corollary
   \ref{inner3}, there exists $\gamma\in\R$ such that
   $$
   e^{\alpha\frac{\vp(z) + 1}{\vp(z) - 1}} = e^{i\gamma} e^{\alpha\frac{z + 1}{z - 1}}.
   $$
   Therefore, we obtain
   $$
   \frac{\vp(z) + 1}{\vp(z) - 1} = i(\gamma + 2k\pi)/\alpha + \frac{z + 1}{z - 1}, k\in\mathbb{Z}.
   $$
   Thus, we have the following expression for $\vp(z)$:
   $$
   	\vp(z)  = \frac{2z + i\zeta (z - 1)}{2 + i\zeta(z - 1)}
    = 1 + \frac{2(z - 1)}{2 + (i\zeta)(z - 1)}
   $$
   where $\zeta = (\gamma + 2k\pi)/\alpha$. Now it follows that $\vp(1) = 1$ and $\vp'(1) = 1$ i.e., $\vp$ is a parabolic disc automorphism
   with fixed point $1$.
   \epf

   From Theorems 6 and 7 of \cite{cowen}, it is worth to note that for $\psi=\exp{(\alpha\frac{z + 1}{z - 1})}, \alpha>0$, it is known that  $\psi H^2$ is invariant under $C_\vp$ if and only if $\vp(1) = 1$ and $\vp'(1) \leq 1$. The same result can be verified for $\psi H^\infty$. As $C_\vp$ is an automorphism $\psi H^\infty$ if and only if $\psi H^\infty$ is invariant under both $C_\vp$ and $C_{\vp^{-1}}$, we get an alternative proof of Theorem \ref{thm58}.
   From Problem 6 of Exercises 0.5 of \cite{Shapiro}, as an immediate consequence of the above theorem, we get the following result.
   \bcor
   Let $\alpha>0$. Then,
   $$\mbox{Automorphisms of } e^{\alpha\frac{z + 1}{z - 1}} H^\infty=\{C_\vp: \vp(z)    = 1 + \frac{2(z - 1)}{2 + (i\zeta)(z - 1)}, \zeta\in \R\}
     $$
   \ecor

   \brem\label{rotat}
  Let $\vp ,\eta\in\mathrm{Aut}(\D)$. Then $C_\vp$ is an automorphism $\psi H^\infty$ if and only if $C_{\eta^{-1}\circ \vp \circ\eta}$ is an automorphism on $\psi \circ \eta H^\infty$.
  \erem
   By Theorem \ref{thm58} and choosing $\eta(z) =\overline{w}z$ in the remark \ref{rotat}, we can deduce the following corollary.
   \bcor
   Fix $|w|=1$ and $\alpha>0$. Then
   \begin{align*}
   \mbox{Automorphisms of } e^{\alpha\frac{z + w}{z - w}} H^\infty &=\{C_{\phi}: \phi(w) = w \mbox{ and } \phi^{'}(w) = 1 \}\\
   &=\{C_{wz\circ\vp\circ\overline{w}z}: \vp(z)= 1 + \frac{2(z - 1)}{2 + (i\zeta)(z - 1)}, \zeta\in \R\}\\
   &=\{C_{\phi}: \phi(z)= w + \frac{2w(z - w)}{2w + (i\zeta)(z - w)}, \zeta\in \R\}.
   \end{align*}
     %
   \ecor

\bprop\label{twoatomsing1min1}
Let $\vp\in \mathrm{Aut}(\D)$ and $\psi(z) = e^{\frac{z + 1}{z - 1}}e^{\frac{z - 1}{z + 1}}=e^{\frac{2(z^2 + 1)}{z^2 - 1}}$. Then $C_\vp$ is an automorphism of $\psi H^\infty$ if and only if $\vp =\pm z$.
\eprop
\bpf
If $\vp=\pm z$, then  it is easy to see that $\psi\circ\vp = \psi$ and hence in each of these cases $C_\vp$ is an automorphism of $\psi H^\infty$.

Conversely, suppose $C_\vp$ is an automorphism of $\psi H^\infty$. As $\psi$ is an inner function, by  Corollary \ref{inner3} we must have
 $\psi\circ\vp = e^{i\gamma}\psi$ for some $\gamma \in\R$. Thus,
 $$
   e^{\frac{2(\vp(z)^2 + 1)}{\vp(z)^2 - 1}} = e^{i\gamma} e^{\frac{2(z^2 + 1)}{z^2 - 1}}.
   $$
   Consequently, we obtain
   $$
   \frac{\vp(z)^2 + 1}{\vp(z)^2 - 1} = i\zeta + \frac{z^2 + 1}{z^2 - 1},
   $$
   where $\zeta=(\gamma+2k\pi)/2$ for some $k\in\mathbb{Z}$. This simplifies to
   \begin{equation}\label{phisqeq}
  \vp(z)^2 = \frac{(2+i\zeta)z^2-i\zeta}{i\zeta z^2 + (2-i\zeta)} 
  \end{equation}

   Now let $\vp(z)=\alpha\frac{a-z}{1-\overline{a}z}$ for $a\in\D,\alpha\in\mathbb{T}$. Substituting this in  \eqref{phisqeq},
    cross multiplying and comparing coefficients of $z^3$ and constant terms on both side, we get
   \begin{equation}\label{zcubecoefs}
   \alpha^2 a (i\zeta) = \overline{a}(2+i\zeta) 
   \end{equation}
   and
   \begin{equation}\label{constcoefs}
   a^2 = \frac{-i\zeta}{\alpha^2(2-i\zeta)}. 
   \end{equation}
   If $a\neq 0$, then by taking the modulus of both side of \eqref{zcubecoefs} we get a contradiction, thus $a=0$. Now from \eqref{constcoefs} we get $\zeta = 0$. Therefore, \eqref{phisqeq} implies $\vp(z)^2 = z^2$ for $z\in\D$. Finally, the use of the identity theorem gives us $\vp=\pm z$.
\epf
By Proposition \ref{twoatomsing1min1} and choosing $\eta(z) =\overline{w}z$ in the remark \ref{rotat}, we can deduce the following corollary.
\bcor
Let $\vp\in \mathrm{Aut}(\D)$ and $\psi(z) = e^{\frac{z + w}{z - w}}e^{\frac{z - w}{z + w}}$, for $w\in\mathbb{T}$. Then $C_\vp$ is an automorphism of $\psi H^\infty$ if and only if $\vp =\pm z$.
\ecor

%
%
%
%
%

    \bprop\label{boundaryzero1} Fix $m\in \mathbb{N}$.
   $$\mbox{Automorphisms of } (z-1)^m H^\infty=\{C_\vp: \vp(z)    = \frac{\overline{a-1}}{a-1}\frac{z-a}{1-\overline{a}z}, a\in\D\}
     $$
   \eprop
   \bpf
   By Theorem \ref{Poly1}, $C_\vp$ is an automorphism of $(z-1)^m H^\infty$ if and only if $\vp(1) = 1$ i.e., $\vp = \frac{1-\overline{a}}{a-1}\tau_a$ for some $a\in\D$.
   \epf

If we choose $\eta(z) =\overline{w}z$ for $w\in\mathbb{T}$ in Remark \ref{rotat}, then from Proposition \ref{boundaryzero1} we can get all automorphisms of $(z-w)^mH^\infty$.

   \bprop\label{boundzero1min1same}
  Fix $m\in \mathbb{N}$. Then $C_\vp$ is an automorphism of $(z^2-1)^m H^\infty$ if and only if $\vp = -\tau_a$ or $\vp = \tau_a$ for some $a \in (-1, 1)$.
   \eprop
   \bpf
   By Theorem \ref{Poly1}, $C_\vp$ is an automorphism $(z^2-1)^m H^\infty$ if and only if any one of the following two cases for $\vp$ occurs:
      \begin{enumerate}
   	\item $\varphi(1) = 1$ and $\varphi(-1) = -1$, that is $\varphi = -\tau_a$ for some $a \in (-1, 1)$.
   	\item $\varphi(1) = -1$ and $\varphi(-1) = 1$, that is $\varphi = \tau_a$ for some $a \in (-1, 1)$.
   \end{enumerate}
   \epf
  The following result is immediate now.

   \bcor\label{boundzero1min1dif}
    $C_\vp$ is an automorphism of $(z-1)^m(z+1)^n H^\infty$, where $m\neq n$, if and only if $\vp = -\tau_a$ for some $a \in (-1, 1)$.
   \ecor


  If we choose $\eta\in\mathrm{Aut}(\D)$ such that $\eta(w_1)=1$ and $\eta(w_2)=-1$ for $w_1,w_2\in\mathbb{T}$ in Remark \ref{rotat}, then from
  Proposition \ref{boundzero1min1same} and Corollary \ref{boundzero1min1dif}, we can get all automorphisms of $(z-w_1)^m(z-w_2)^mH^\infty$ and $(z-w_1)^m(z-w_2)^nH^\infty$, respectively.


  Now, we give another simpler proof of Proposition \ref{coroH_an} for the case of $\mathcal{A}_a^n$.

   \bprop\label{cor2}
    $C_\vp$ is an automorphism $\tau_{a}^nH^{\infty}$ if and only if there exists $\theta \in \mathbb{R}$ such that $\vp =\tau_{a}\circ e^{i\theta}z\circ\tau_{a}$.
   \eprop

   \bpf
    By Corollary \ref{anyblaschke}, $C_\vp$ is an automorphism if and only if  $\vp(a) = a$. Since $\vp$ is disc automorphism, it must be of the form $ \vp = \tau_{a}\circ e^{i\theta}z\circ \tau_{a}$ for some $\theta \in \mathbb{R}$.
   \epf

 We are now ready for the characterization of automorphism of subalgebra\\ $(\tau_a)^m(\tau_b)^nH^\infty$, where $a \neq b \in \D, m, n\in \mathbb{N}$.

   \bprop \label{cor4}
   Let  $a,b\in\D$ be distinct points. Then $C_\vp$ is an automorphism of $(\tau_a\tau_b)^mH^\infty$ if and only if either $\vp$ is identity or $\vp = \tau_a \circ \tau_c \circ \tau_a$, with  $ c = \tau_a(b)$.
   \eprop

   \bpf
   By Corollary \ref{anyblaschke}, $C_\vp$ is an automorphism of $(\tau_a\tau_b)^mH^\infty$ if and only if $\vp(a) = a$ and $\vp(b) = b$ or $ \vp(a) = b$ and $\vp(b) = a$.
In the former case,    $\vp$ must be the identity map on $\D$. Now, suppose
    $ \vp(a) = b$ and $\vp(b) = a$. We define $\phi = \tau_a \circ \tau_c \circ \tau_a\circ \vp$, where $c = \tau_{a}(b)$. This implies that $\phi$ is a disc automorphism with two fixed points $a$ and $b$ in $\D$. Therefore $\phi$ must be the identity map on $\D$, and hence $\vp = \tau_a \circ \tau_c \circ \tau_a$. This completes the proof.
   \epf

   \bcor \label{cor5}
   Let $a,b\in\D$ be distinct and $m\neq n$. Then $C_\vp$ is an automorphism of $(\tau_a)^m(\tau_b)^n H^{\infty}$ if and only if $\vp$ is identity.
   \ecor

  \textbf{Concluding Remarks}: In Section \ref{chara auto}, we proved that every automorphism of $\psi H^\infty$  or $\psi A(\D)$ is a composition operator
  induced by a disc automorphism. 
  We also proved in Section \ref{sectionj} that every automorphism of $\mathcal{B}_a^n$ or $\mathcal{\tilde{B}}_a^n$, which are not of the form $\psi H^\infty$, is also a composition operator induced by a disc automorphism.
  It leads us to expect  that every automorphism of any subalgebra of $H^\infty$  can be extended to an automorphism of $H^\infty$, and hence it will be composition operator. In view of this, we make the following conjecture.

  \textbf{Conjecture}: Let $\mathcal{A}$ be any subalgebra of bounded analytic functions $H^\infty$ on the unit disc in the complex plane. Every algebra automorphism of $\mathcal{A}$ must be a composition operator induced by some disc automorphism.\\
{\bf Acknowledgments.}
The authors thank the referee for many useful comments which improve the
presentation considerably.
\bibliographystyle{plain}
 \bibliography{allrefs}

\end{document}